\documentclass[reqno,a4paper,10pt]{amsart}
\usepackage{enumitem}
\usepackage{amsmath, amssymb, amsthm}
\usepackage[utf8]{inputenc}
\usepackage{url}
\usepackage[english]{babel}

\thanks{An extended
	abstract of this paper will be published in the Proceedings of the 31st International Conference on Probabilistic, Combinatorial and Asymptotic Methods for the Analysis of Algorithms (AofA 2020), pp. 4:1-4:17.}
\author{Mihyun Kang, Michael Missethan}
\address{Graz University of Technology, Institute of Discrete Mathematics, Steyrergasse 30, 8010 Graz, Austria}
\email{\{kang,missethan\}@math.tugraz.at}
\thanks{Supported by Austrian Science Fund (FWF): I3747 and W1230}

\title[Sparse random outerplanar graphs]{The giant component and 2-core in sparse random outerplanar graphs}


\newtheorem{thm}{Theorem}[section]

\newtheorem{lem}[thm]{Lemma}
\newtheorem{conjecture}[thm]{Conjecture}

\theoremstyle{definition}
\newtheorem{definition}[thm]{Definition}
\begin{document}

\begin{abstract}
Let $A(n,m)$ be a graph chosen uniformly at random from the class of all vertex-labelled outerplanar graphs with $n$ vertices and $m$ edges. We consider $A(n,m)$ in the sparse regime when $m=n/2+s$ for $s=o(n)$. We show that with high probability the giant component in $A(n,m)$ emerges at $m=n/2+O\left(n^{2/3}\right)$ and determine the typical order of the 2-core. In addition, we prove that if $s=\omega\left(n^{2/3}\right)$, with high probability every edge in $A(n,m)$ belongs to at most one cycle.
\end{abstract}

\maketitle

\section{Introduction}

\subsection{Motivation}
In 1959 Erd{\H o}s and R\'{e}nyi \cite{erdoes1} introduced the so-called \textit{Erd{\H o}s-R\'{e}nyi graph} $G(n,m)$, a graph chosen uniformly at random from the class of all vertex-labelled graphs on vertex set $\{1, \ldots, n\}$ with $m=m(n)$ edges. Since then, the asymptotic behaviour of $G(n,m)$ was extensively studied (see e.g. \cite{rg3, rg1, rg2}). In particular, it was investigated how the component structure of $G(n,m)$ changes, when $m=m(n)$ varies and whether there are ranges of $m$, where this change is very significant. Such dramatic changes are called \textit{phase transitions}. For example, Erd{\H o}s and R\'{e}nyi \cite{erdoes2} showed that the order (that is, the number of vertices) of the largest component in $G(n,m)$ changes drastically when $m \sim n/2$. Later Bollobás \cite{general1} and {\L}uczak \cite{general2} looked more closely at the critical range $m=n/2+o(n)$. 

Throughout the paper, we denote the components of a graph $G$ by $H_1=H_1(G), H_2=H_2(G), \ldots$ in such a way that $|H_i|\geq |H_j|$, whenever $i \leq j$, where $|H_i|$ is the number of vertices in $H_i$. In addition, we use the asymptotic notation from \cite{asy_not}.

\begin{thm}[\cite{general1, general2}]\label{thm:general}
Let $m=n/2+s$, where $s=s(n)=o(n)$ and let $G=G(n,m)$. Then for every $i \in \mathbb{N}$ the following holds with high probability\footnote{With probability tending to 1 as $n$ tends to infinity, whp for short.}.
\begin{enumerate}[label=(\roman*)]
\item 
If $\frac{s^3}{n^2}\to -\infty$, then $H_i$ is a tree and $\left|H_i\right|=\left(1/2+o(1)\right)\frac{n^2}{s^2}\log \frac{|s|^3}{n^2}.$
\item
If $\frac{s^3}{n^2}\to c \in \mathbb{R}$, then $|H_i|=\Theta_p\left(n^{2/3}\right).$
\item
If $\frac{s^3}{n^2}\to \infty$, then $\left|H_1\right|=\left(4+o(1)\right)s.$ For $i\geq 2$, we have $|H_i|=o\left(n^{2/3}\right)$.
\end{enumerate}
\end{thm}

This drastic change of the component structure at $m=n/2+O\left(n^{2/3}\right)$ is called the \emph{emergence of the giant component}. These results raised the question whether there are also phase transitions in other classes of random graphs. {\L}uczak and Pittel \cite{luczak_pittel_1992} considered this question for $F(n,m)$, a graph chosen uniformly at random from all vertex-labelled forests with $n$ vertices and $m$ edges. They showed that, analogous to $G(n,m)$, the giant component in $F(n,m)$ emerges at $m=n/2+O\left(n^{2/3}\right)$. Kang and {\L}uczak \cite{planar} showed that the same is true for $P(n,m)$, a graph chosen uniformly at random from all vertex-labelled planar graphs with $n$ vertices and $m=m(n)$ edges. Later Kang, Mo{\ss}hammer, and Spr\"ussel \cite{surface} extended this result even  to graphs on orientable surfaces. 

Surprisingly, this problem for a random \emph{outerplanar} graph is still open, although the class of outerplanar graphs lies \lq between\rq\ the class of forests and the class of planar graphs and therefore we expect similar behaviours. (A graph is outerplanar if it has an embedding in the plane in such a way that every vertex lies on the outer face, equivalently, a graph is outerplanar iff it  contains neither $K_4$ nor $K_{2,3}$ as a minor.)  In this paper we solve this open problem on the emergence of the giant component in a random outerplanar graph.  

Kang, Mo{\ss}hammer, and Spr\"ussel \cite{surface} used the core-kernel approach to obtain their results on the giant component in $S_g(n,m)$, a graph chosen uniformly at random from all vertex-labelled graphs with $n$ vertices,  $m=m(n)$ edges and genus at most $g$ (for any constant $g\ge 0$). This method is mainly based on the following decomposition. We call a component of a graph $G$ \emph{complex} if it has at least two cycles. We decompose $G$ into the \emph{complex part} $Q_G$, which is the union of all complex components, and into non-complex components. Then we extract the \emph{core} $C_G$, which is the maximal subgraph of $Q_G$ of minimum degree at least two. Finally, we consider the \emph{kernel} $K_G$, which can be obtained from $C_G$ by the following operation. Every maximal path $P$ consisting of vertices of degree two is replaced by an edge between the vertices of degree at least three that are adjacent to the end vertices of $P$. Conversely, starting from kernels (as base cases) we can construct cores by subdividing edges with additional vertices. Similarly, the complex part can be formed by replacing every vertex in the core by a rooted tree. Finally, we obtain the whole graph $G$ by choosing the complex part and non-complex components. 

However, we \emph{cannot} apply the core-kernel approach to {\em outerplanar} graphs, because this method is mainly based on the fact that a graph $G$ is embeddable on a surface if and only if its kernel $K_G$ is. But an analogous statement for outerplanar graphs is not true, since a subdivision of an outerplanar graph is not necessarily  outerplanar. Therefore, in this paper we shall start directly from cores (as base cases), not from the kernels. One of key steps in this direct core approach is to investigate how the number of outerplanar cores (and complex parts, respectively) changes by addition of a vertex and an edge.  Using our core approach we prove that the giant component in a random outerplanar graph with  $n$ vertices and $m=m(n)$ edges emerges at $m=n/2+O\left(n^{2/3}\right)$. 

\subsection{Main results}
To state our main results we need to introduce some notations. 
Given a graph $G$, we define the excess of a complex component of $G$ to be the difference between the number of its edges and the number of its vertices. The excess of $G$, denoted by $ex(G)$ or $\ell(G)$, is the sum of the excesses of all complex components of $G$. In addition, we denote by $n_C(G)$ the number of vertices in the core $C_G$. 
Let $A(n,m)$ denote a graph chosen uniformly at random from all vertex-labelled outerplanar graphs with $n$ vertices and $m=m(n)$ edges. 

\begin{thm}\label{thm:pt_out}
Let $m=n/2+s$, where $s=s(n)=o(n)$ and let $G=A(n,m)$. For every $i \in \mathbb{N}$ whp the following holds.
\begin{enumerate}[label=(\roman*)]
\item\label{thm:main_1_1}
If $\frac{s^3}{n^2}\to -\infty$, then $H_i$ is a tree and $\left|H_i\right|=\left(1/2+o(1)\right)\frac{n^2}{s^2}\log \frac{|s|^3}{n^2}.$
\item\label{thm:main_1_2} 
If $\frac{s^3}{n^2}\to c \in \mathbb{R}$, then $|H_i|=\Theta_p\left(n^{2/3}\right).$
\item\label{thm:main_1_3}
If $\frac{s^3}{n^2}\to \infty$, then $\left|H_1\right|=2s+O_p\left(n^{2/3}\right)$. For $i\geq 2$, we have $|H_i|=\Theta_p \left(n^{2/3}\right)$. We also have $n_C(G)=\Theta\left(sn^{-1/3}\right)$ and $ex(G)=\Theta\left(sn^{-2/3}\right)$.
\end{enumerate}
\end{thm}

To prove Theorem \ref{thm:pt_out} we shall use some auxiliary results about {\em cactus} graphs, which form a subfamily of the class of outerplanar graphs and are interesting in their own -- a cactus graph is a graph in which every edge belongs to at most one cycle. A simple, but important observation is that a graph is a cactus graph if and only if its kernel is a cactus graph. Therefore, analogously to the case of random graphs on surfaces \cite{surface}  we can apply the aforementioned core-kernel approach to obtain results on the component structure of a random cactus graph, such as the order of the largest component, the core, and the kernel. In addition, we determine the asymptotic number of cubic (i.e. 3-regular) cactus multigraphs using singularity analysis of generating functions which arise from the standard decomposition of graphs into smaller building blocks.

We denote by $T(n,m)$ a graph chosen uniformly at random from all vertex-labelled cactus graphs with $n$ vertices and $m=m(n)$ edges. In addition, let $\mathcal{K}(2n,3n)$ be the class of all cubic cactus {\em weighted multigraphs} with $2n$ vertices and $3n$ edges, and $\mathcal{K}_c(2n,3n)$ be the subclass of $\mathcal{K}(2n,3n)$ containing all connected graphs. Here every multigraph $K$ is counted with a weight of $w(K)=2^{-e_1(K)-e_2(K)}$, where $e_1(K)$ denotes the number of loops in $K$ and $e_2(K)$ the number of double edges (see \cite[p.5]{comp_fac} for details of the weight of a multigraph).

\begin{thm}\label{thm:pt_cac}

\begin{enumerate}[label=(\roman*)]
\item\label{thm:pt_cac_1}
Let $m=n/2+s$, where $s=s(n)$, $n^{2/3}\ll s \ll n$ and $G=T(n,m)$. Then whp $\left|H_1\right|=2s+O_p\left(n^{2/3}\right)$, $n_C(G)=\Theta\left(sn^{-1/3}\right)$, $ex(G)=\Theta\left(sn^{-2/3}\right)$, and the kernel $K_G$ is cubic. \item\label{thm:pt_cac_2}
There are constants $c_0, c_1, \gamma>0$ such that as $n\to \infty$,
\begin{align*}
|\mathcal{K}(2n,3n)|&=(1+o(1))c_0n^{-5/2} \gamma^{2n}(2n)!, \\
\text{and }~~~~~~~~~~~|\mathcal{K}_c(2n,3n)|&=(1+o(1))c_1n^{-5/2} \gamma^{2n}(2n)!.~~~~~~~~~~~
\end{align*}
\end{enumerate}
\end{thm}

Finally, we use Theorem \ref{thm:pt_out} to show that when $m=n/2+s$ for  $n^{2/3}\ll s \ll n$, the two random graphs  $A(n,m)$ and $T(n,m)$ are \lq contiguous\rq, meaning that they are indistinguishable in view of properties that hold whp. Such a contiguity of two models will turn out to be very helpful for further investigations of the behaviour of $A(n,m)$, partly because the core-kernel approach is applicable for $T(n,m)$.

\begin{thm}\label{thm:out_cac}
Let $m=n/2+s$, where $s=s(n)$ and $n^{2/3}\ll s \ll n$. Then, whp every edge in $A(n,m)$ belongs to at most one cycle. In other words, whp  $A(n,m)$ is a cactus graph.
\end{thm}

\section{Proof strategy of Theorem \ref{thm:pt_out}}\label{sec:proof_strategy}
We start with the cases $s^3/n^2 \to -\infty$ and $s^3/n^2\to c\in \mathbb{R}$. By a well-known fact (see Lemma \ref{lem:unicyclic}\ref{lem:unicyclic_1},\ref{lem:unicyclic_2}) we obtain
$
\liminf_{n \to \infty} \mathbb{P}\left[{G(n,m) \text{ is outerplanar}}\right]>0.
$
Thus, each property that holds whp in $G(n,m)$ is also true whp in $A(n,m)$ and the Statements \ref{thm:main_1_1} and \ref{thm:main_1_2} follow from Theorem \ref{thm:general}. Thus, it suffices to prove \ref{thm:main_1_3}, for which we use the direct core approach. To illustrate this approach, we introduce further notations.  
\begin{definition}
We denote by
\begin{itemize}
\item
$\mathcal{A}$ the class of all outerplanar graphs;
\item
$\mathcal{Q}$ the class of all complex outerplanar graphs (i.e. complex parts of graphs in $\mathcal{A}$);
\item
$\mathcal{C}$ the class of all complex outerplanar graphs with minimum degree at least two (i.e. cores of graphs in $\mathcal{A}$);
\item
$\mathcal{U}$ the class of all graphs without complex components.
\end{itemize}
In addition, for any graph class $\mathcal{X}$ we denote by $\mathcal{X}(n,m)$ the subclass containing those graphs with $n$ vertices and $m$ edges.
\end{definition}

\begin{definition}
Let $G$ be a graph with $n$ vertices and $m$ edges. We denote by
\begin{itemize}
\item
$n_Q=n_Q(G)$ the number of vertices in the complex part $Q_G$;
\item
$n_C=n_C(G)$ the number of vertices in the core $C_G$;
\item 
$\ell=\ell(G)$ the excess of $G$, i.e. the difference between the number of edges and the number of vertices in the complex part $Q_G$;
\item
$n_U=n_U(G):=n-n_Q$ the number of vertices in $G$ outside the complex part $Q_G$;
\item
$m_U=m_U(G):=m-n_Q-\ell$ the number of edges in $G$ outside the complex part $Q_G$ (with $n_Q$ vertices and $n_Q+\ell$ edges).
\end{itemize}
\end{definition} 

We reverse the decomposition in the core approach to obtain relations between the classes defined above. We observe that each outerplanar graph can be constructed in a unique way by combining a complex graph and non-complex components. Similarly, a complex graph can be formed by choosing the core and replacing each vertex of the core by a rooted tree. It is well known that we have $n_Cn_Q^{n_Q-n_C-1}$ different possibilities for choosing these trees (see e.g. \cite{TAKACS}). Hence, we obtain
\begin{align}
\left|\mathcal{A}(n,m)\right|
&=\sum\nolimits_{n_Q, \ell}\binom{n}{n_Q}|\mathcal{Q}(n_Q, n_Q+\ell)|\cdot|\mathcal{U}(n_U, m_U)|
=\sum\nolimits_{n_Q, \ell}\tau(n_Q,\ell), \label{eq:1}\\
|\mathcal{Q}(n_Q,n_Q+\ell)|
&=\sum\nolimits_{n_C}\binom{n_Q}{n_C}|\mathcal{C}(n_C, n_C+\ell)|n_Cn_Q^{n_Q-n_C-1}
=\sum\nolimits_{n_C}\rho(n_C), \label{eq:2}
\end{align}
where we define 
\begin{align*}
 \tau(n_Q,\ell)&:=\binom{n}{n_Q}|\mathcal{Q}(n_Q, n_Q+\ell)|\cdot|\mathcal{U}(n_U, m_U)|,\\
 \rho(n_C)&:=\binom{n_Q}{n_C}|\mathcal{C}(n_C, n_C+\ell)|n_Cn_Q^{n_Q-n_C-1}. 
\end{align*}
In the sums of (\ref{eq:1}) and (\ref{eq:2}) we did not specify precisely in which sets the summation indices lie. But it is convenient to consider only terms, which are non-zero. We call the corresponding indices \emph{admissible}. The next step is to find in the sums (\ref{eq:1}) and (\ref{eq:2}) those terms, which are significantly larger than the other ones. In order to make that more precise, we use the following terminology.

\begin{definition}\label{def:gen_sum}
For each $n \in \mathbb{N}$ let $I_0(n), I(n) \subseteq \mathbb{N}$ be finite index sets such that $I_0(n) \subseteq I(n)$. In addition, let $\sigma_n(i)\geq0$ for each $i \in I(n)$. Then the \textit{main contribution} to the sum $\sum_{i\in I(n)}\sigma_n(i)$ is provided by $i \in I_0(n)$ if $\sum_{i \in I(n)\backslash I_0(n)}\sigma_n(i)=o\left(\sum_{i \in I(n)}\sigma_n(i)\right)$
for $n \to \infty$. In that case, we also say that the terms provided by $i \in I(n)\backslash I_0(n)$ are \textit{negligible}.
\end{definition}

Now the goal is to find sets $I_{n_Q}, I_\ell$ and $I_{n_C}$ such that the main contributions to \eqref{eq:1} and \eqref{eq:2} are provided by $n_Q \in I_{n_Q}, \ell\in I_\ell$, and  $n_C\in I_{n_C}$. Having such sets we immediately get results about the structure of a random outerplanar graph $G=A(n,m)$. Namely, that whp $n_Q(G) \in I_{n_Q}, \ell(G)\in I_\ell$, and $n_C(G)\in I_{n_C}$. To get strong results, we aim to find sets $I_{n_Q}, I_\ell$, and $I_{n_C}$, which are as small as possible. Afterwards we use this concentration information and a double counting argument (see Lemma \ref{lem:giant_complex}) to deduce the component structure of $G$. The main challenge  is to determine $I_{n_Q}, I_\ell$, and $I_{n_C}$.

In order to illustrate our main idea of the analysis of the sums \eqref{eq:1} and \eqref{eq:2}, we consider the generic sums 
$\Sigma_n=\sum_{i\in I(n)}\sigma_n(i)$ from Definition \ref{def:gen_sum}. The goal is to find \lq small\rq~sets $I_0(n)$ such that the main contribution to $\Sigma_n$ is provided by $i\in I_0(n)$ or equivalently \lq large\rq~sets $I_1(n)$ such that the terms provided by $i\in I_1(n)$ are negligible in $\Sigma_n$. Our method to find these sets $I_1(n)$ is mainly based on the following observation. 

\begin{lem}\label{lem:main_analysis}
For each $n \in \mathbb{N}$ let $I_1(n), I(n) \subseteq \mathbb{N}$ be finite index sets such that $I_1(n) \subseteq I(n)$ and let $\sigma_n(i)\geq0$ for each $i\in I(n)$. In addition, for each $n \in \mathbb{N}$ let $f_n: I_1(n)\to I(n)$ be a function. We assume that there are a function $\varepsilon$ with $\varepsilon(n)=o(1)$ and a constant $M>0$ such that for all $n \in \mathbb{N}, i \in I_1(n)$ and $j\in I(n)$
\begin{align}
\frac{\sigma_n(i)}{\sigma_n\left(f_n(i)\right)}&\leq \varepsilon(n), \label{lem:main_analysis_1} \\
\text{and}~~~~\left|f_n^{-1}\left(\left\{j\right\}\right)\right|&\leq M.~~~~\label{lem:main_analysis_2}
\end{align}
Then the terms provided by $i \in I_1(n)$ are negligible in $\sum_{i\in I(n)}\sigma_n(i)$.
\end{lem}

In most cases when we apply Lemma \ref{lem:main_analysis}, the functions $f_n$ will be of the form $f_n(i)=i+g(n)$ for some function $g:\mathbb{N}\to \mathbb{Z}$ or of the form $f_n(i)=\lfloor \delta i \rfloor$ for some constant $\delta>0$. We note that such functions $f_n$ always fulfil (\ref{lem:main_analysis_2}) for some $M>0$. Thus, it remains to find a function $\varepsilon$ with $\varepsilon(n)=o(1)$ such that (\ref{lem:main_analysis_1}) is satisfied. For simplicity, we demonstrate our method of doing that only for the case when $f_n(i)=i+g(n)$ for some function $g$ with $g(n)>0$. Moreover, we assume that $I(n)=\{a_n, a_n+1, \ldots, b_n\}$ for some $a_n<b_n$. We observe that
\begin{align}\label{lem:main_analysis3}
\frac{\sigma_n(i)}{\sigma_n(f_n(i))}=\frac{\sigma_n(i)}{\sigma_n(i+g(n))}=\prod_{k=i}^{i+g(n)-1}\frac{\sigma_n(k)}{\sigma_n(k+1)}.
\end{align}
Thus, we aim to find good upper bounds for $\frac{\sigma_n(k)}{\sigma_n(k+1)}$. We commonly state these bounds in the form $\exp(h(n))$ for some function $h:\mathbb{N}\to\mathbb{R}$. Then, if we assume 
\begin{align}\label{lem:main_analysis4}
\frac{\sigma_n(k)}{\sigma_n(k+1)}\leq \exp(h(n)), ~~~~\forall n \in \mathbb{N}, \forall k \in \left\{i, \ldots, i+g(n)-1\right\},
\end{align}
we get together with (\ref{lem:main_analysis3}), $\frac{\sigma_n(i)}{\sigma_n(f_n(i))}\leq \exp\left(g(n)h(n)\right).$ If we find such functions $g$ and $h$ with $g(n)h(n) \to -\infty$ for $n\to \infty$, then we can apply Lemma \ref{lem:main_analysis} (see Appendix \ref{sec:appendix} for an application of Lemma \ref{lem:main_analysis}). We can summarise the above idea as follows. The key for a good analysis of the sum $\sum_{i\in I(n)}\sigma_n(i)$ is to have good bounds for the fractions $\frac{\sigma_n(k)}{\sigma_n(k+1)}$ or equivalently good bounds for $\frac{\sigma_n(k+1)}{\sigma_n(k)}$. 

Now we describe how we find these bounds for the sums in \eqref{eq:1} and \eqref{eq:2}. 
In order to find good bounds for $\frac{\rho(n_C+1)}{\rho(n_C)}$, it suffices to estimate $\frac{\left|\mathcal{C}\left(n_C+1, n_C+1+\ell\right)\right|}{\left|\mathcal{C}(n_C,n_C+\ell)\right|}$   (see Lemma \ref{lem:estimate_core}). To that end, we construct graphs in $\mathcal{C}\left(n_C+1, n_C+1+\ell\right)$ as follows: Let $H \in \mathcal{C}(n_C,n_C+\ell)$ and an edge $e$ of $H$ be given. Then we obtain in \lq most\rq\ cases a graph $H' \in \mathcal{C}(n_C+1,n_C+1+\ell)$ if we subdivide $e$ by one vertex and label this new vertex with $n_C+1$. By a careful analysis of this construction we will obtain good estimates for $\frac{\rho(n_C+1)}{\rho(n_C)}$. 

In the next step we consider the sum in (\ref{eq:1}) and shall determine $I_{n_Q}$ and $I_\ell$. To that end, we look at the fractions $\frac{\tau(n_Q+1,\ell)}{\tau(n_Q,\ell)}$ and $\frac{\tau\left(n_Q,\left \lfloor \delta\ell\right \rfloor\right)}{\tau(n_Q,\ell)}$ for a constant $\delta>0$. To get bounds for the term $|\mathcal{U}(n_U, m_U)|$, we will use Lemma \ref{lem:unicyclic}. Thus, it remains to find estimates for $\frac{|\mathcal{Q}(n_Q+1, n_Q+1+\ell)|}{|\mathcal{Q}(n_Q, n_Q+\ell)|}$ and $\frac{|\mathcal{Q}(n_Q, n_Q+\left \lfloor \delta\ell\right \rfloor)|}{|\mathcal{Q}(n_Q, n_Q+\ell)|}$. 
For the first fraction (see Lemma \ref{lem:bound_complex}) we define for $i \in \{0, 1\}$
\[
\rho_i(n_C)=\rho_i(n_C, n_Q, \ell):=\binom{n_Q+i}{n_C}\left|\mathcal{C}(n_C, n_C+\ell)\right|n_C(n_Q+i)^{n_Q+i-n_C-1}.
\]
With this notation we have  
\begin{align}\label{eq:frac_comp}
\frac{|\mathcal{Q}(n_Q+1, n_Q+1+\ell)|}{|\mathcal{Q}(n_Q, n_Q+\ell)|}=\frac{\sum_{n_C}\rho_1(n_C)}{\sum_{n_C}\rho_0(n_C)}.
\end{align}
From the analysis of (\ref{eq:2}) we already know sets $I_0, I_1$ such that the main contributions to $\sum_{n_C}\rho_0(n_C)$ and $\sum_{n_C}\rho_1(n_C)$ are provided by $n_C \in I_0$ and $n_C \in I_1$, respectively. We will see that we may assume $I:=I_0=I_1$. Then we will get a good bound for (\ref{eq:frac_comp}) if for $n_C \in I$ we estimate the fraction
\begin{align}\label{eq:frac_r}
\frac{\rho_1(n_C)}{\rho_0(n_C)}=\frac{\left(n_Q+1\right)^2}{n_Q-n_C+1}\left(\frac{n_Q+1}{n_Q}\right)^{n_Q-n_C-1}.
\end{align}	
For the fraction $\frac{|\mathcal{Q}(n_Q, n_Q+\left \lfloor \delta\ell\right \rfloor)|}{|\mathcal{Q}(n_Q, n_Q+\ell)|}$ (see Lemma \ref{lem:frac_l}), we will use that
\begin{align}\label{eq:bound_complex}
\frac{|\mathcal{Q}_C(n_Q, n_Q+\left \lfloor \delta\ell\right \rfloor)|}{|\mathcal{Q}_P(n_Q, n_Q+\ell)|}\leq\frac{|\mathcal{Q}(n_Q, n_Q+\left \lfloor \delta\ell\right \rfloor)|}{|\mathcal{Q}(n_Q, n_Q+\ell)|}\leq \frac{|\mathcal{Q}_P(n_Q, n_Q+\left \lfloor \delta\ell\right \rfloor)|}{|\mathcal{Q}_C(n_Q, n_Q+\ell)|},
\end{align}
where $\mathcal{Q}_P(n_Q, n_Q+\ell)$ denotes the class of all complex planar graphs with $n_Q$ vertices and $n_Q+\ell$ edges and $\mathcal{Q}_C(n_Q, n_Q+\ell)$ the class of all complex cactus graphs with $n_Q$ vertices and $n_Q+\ell$ edges. We get estimates for $\left|\mathcal{Q}_C(n_Q, n_Q+\ell)\right|$ and $\left|\mathcal{Q}_P(n_Q, n_Q+\ell)\right|$ by using the core-kernel approach (see Lemmas \ref{lem:complex_cacti} and \ref{lem:complex_planar}). In order to show that these bounds are tight enough, we make the following observations. We will see that there is a constant $c>0$ such that
\begin{align}\label{error}
\frac{\left|\mathcal{Q}_P(n_Q, n_Q+\ell)\right|}{\left|\mathcal{Q}_C(n_Q, n_Q+\ell)\right|}\leq c^\ell,
\end{align}
Thus, we make a multiplicative error of at most $c^\ell$ if we use $\left|\mathcal{Q}_P(n_Q, n_Q+\ell)\right|$ as an estimate for $\left|\mathcal{Q}(n_Q, n_Q+\ell)\right|$. We observe that the possible error increases at most by the constant factor $c$ if we increase $\ell$ by one. On the other hand, we will get  
$\frac{\tau(n_Q,\ell+1)}{\tau(n_Q,\ell)}\approx \Theta(1)\frac{n_Q^{3/2}}{\ell^{3/2}}\frac{1}{n}.$  
Hence, $\tau(n_Q,\ell)$ decays in $\ell$ outside the range $\ell=\Theta\left(n_Qn^{-2/3}\right)$ \lq much faster\rq~than the growth of the error in (\ref{error}).
Having found estimates for $\frac{|\mathcal{Q}(n_Q+1, n_Q+1+\ell)|}{|\mathcal{Q}(n_Q, n_Q+\ell)|}$ and $\frac{|\mathcal{Q}(n_Q, n_Q+\left \lfloor \delta\ell\right \rfloor)|}{|\mathcal{Q}(n_Q, n_Q+\ell)|}$, we obtain bounds for $\frac{\tau(n_Q+1,\ell)}{\tau(n_Q,\ell)}$ and $\frac{\tau\left(n_Q,\left \lfloor \delta\ell\right \rfloor\right)}{\tau(n_Q,\ell)}$. Then we can apply Lemma \ref{lem:main_analysis} to find $I_{n_Q}$ and $I_\ell$.

\section{Cores and complex parts: proof of Theorem \ref{thm:pt_out}}
We recall that  for a given graph $G$  we denote by $n_C$ the number of vertices in the core $C_G$ and by $\ell$ the excess of $G$. In addition, $\mathcal{C}$ is the class of all outerplanar cores. Now we use the ideas presented in Section \ref{sec:proof_strategy} and start by finding $I_{n_C}$. To that end, we obtain the following estimates for $\frac{|\mathcal{C}(n_C+1,n_C+1+\ell)|}{|\mathcal{C}(n_C,n_C+\ell)|}$.

\begin{lem}\label{lem:estimate_core}
\begin{enumerate}[label=(\roman*)]
\item\label{lem:estimate_core_1}
For all admissible $n_C$ and $\ell$ we have
\[
\frac{|\mathcal{C}(n_C+1,n_C+1+\ell)|}{|\mathcal{C}(n_C,n_C+\ell)|}\geq n_C+\frac{\ell}{80}.
\]
\item\label{lem:estimate_core_2}
If in addition $n_C-8\ell\geq 0$, then
\[
\frac{|\mathcal{C}(n_C+1,n_C+1+\ell)|}{|\mathcal{C}(n_C,n_C+\ell)|}\leq (n_C+\ell)\frac{n_C+1}{n_C+1-8\ell}.
\]
\end{enumerate}
\end{lem}

Using Lemma \ref{lem:estimate_core} we obtain bounds for $\frac{\rho(n_C+1)}{\rho(n_C)}$, which we can use to analyse the sum in (\ref{eq:2}) and find $I_{n_C}$. The following two lemmas state that we can choose $I_{n_C}=\Theta\left(\sqrt{n_Q\ell}\right)$, provided that $\ell=\omega(1)$. In  Lemmas \ref{lem:sum_1} and \ref{lem:sum_2} we shall see that we may assume $\ell=\omega(1)$.

\begin{lem}\label{lem:bound_r1}
There are $b, c>0$ such that for all admissible $n_Q$ and $\ell$, we have
\[
\sum\nolimits_{n_C\leq c\sqrt{n_Q\ell}} \rho(n_C) \leq \exp(-b \ell)\sum\nolimits_{n_C} \rho(n_C).
\]
\end{lem}

\begin{lem}\label{lem:bound_r2}
For all admissible $n_Q, \ell$ and $c \geq 14$, we have
\[
\sum\nolimits_{n_C\geq c\sqrt{n_Q\ell}}\rho(n_C)\leq \exp \left(-\frac{c}{2}\ell\right)\sum\nolimits_{n_C}\rho(n_C).
\]
\end{lem}

Next, we recall that $\mathcal{U}$ is the class of all graphs without complex components and $\mathcal{Q}$ the class of all complex outerplanar graphs. In addition, for a given graph $G$ we denote by $n_Q$ the number of vertices in the complex part $Q_G$, by $n_U$ the number of vertices outside the complex part and by $m_U$ the number of edges outside the complex part. We aim to find $I_{n_Q}$ and $I_{\ell}$ by analysing $\frac{\tau(n_Q+1,\ell)}{\tau(n_Q,\ell)}$ and $\frac{\tau\left(n_Q,\left \lfloor \delta\ell\right \rfloor\right)}{\tau(n_Q,\ell)}$. To that end, we need the following estimates for $|\mathcal{U}(n_U,m_U)|$.

\begin{lem}[\cite{uni, comp_fac, surface}]\label{lem:unicyclic}
Let $m=n/2+s$ with $s=s(n)<n/2$ and $u(n,m):=|\mathcal{U}(n, m)|{\binom{\binom{n}{2}}{m}}^{-1}$. Then there is a constant $c>0$ such that for
\begin{align*}
f(n,m):=c\left(\frac{2}{e}\right)^{2m-n}\frac{m^{m+1/2}n^{n-2m+1/2}}{\left(n-m\right)^{n-m+1/2}},
\end{align*}
we have
\begin{enumerate}[label=(\roman*)]
\item\label{lem:unicyclic_1}
$u(n,m) \to 1$, if $\frac{s^3}{n^2} \to -\infty$;
\item\label{lem:unicyclic_2}
for each $a \in \mathbb{R}$, there exists a constant $b>0$ such that $u(n,m)\geq b$, whenever $s\leq an^{2/3};$
\item\label{lem:unicyclic_3}
$u(n,m)\leq n^{-1/2}f(n,m)$, if $0<s\leq\frac{n^{3/4}}{2}$;
\item\label{lem:unicyclic_4}
$u(n,m)\leq f(n,m)$, if $s>0$.
\end{enumerate}
\end{lem}

In addition, we use Lemmas \ref{lem:bound_r1} and \ref{lem:bound_r2} and equation (\ref{eq:frac_r}) to obtain estimates for $\frac{|\mathcal{Q}(n_Q+1, n_Q+1+\ell)|}{|\mathcal{Q}(n_Q, n_Q+\ell)|}$.

\begin{lem}\label{lem:bound_complex}
There exist constants $a_1, a_2, \varepsilon>0$ and $K\in \mathbb{N}$ such that for all admissible $n_Q$ and $\ell$ with $K\leq \ell \leq \varepsilon n_Q$, we have
\begin{align*}
\left(n_Q+1\right)\exp\left(1+a_1\frac{\ell}{n_Q}\right)&\leq \frac{|\mathcal{Q}(n_Q+1, n_Q+1+\ell)|}{|\mathcal{Q}(n_Q, n_Q+\ell)|}
\\
&\leq \left(n_Q+1\right)\exp\left(1+a_2\frac{\ell}{n_Q}\right).
\end{align*}
\end{lem}

Next, we estimate $\frac{|\mathcal{Q}(n_Q, n_Q+\left \lfloor \delta\ell\right \rfloor)|}{|\mathcal{Q}(n_Q, n_Q+\ell)|}$ by using (\ref{eq:bound_complex}). To that end, we need the following two results, which can be obtained by using the core-kernel approach.

\begin{lem}\label{lem:complex_cacti}
There exist constants  $a_1, a_2, \gamma, K, \varepsilon>0$ and $b_1, b_2 \in \mathbb{R}$ such that for all admissible $n_Q$ and $\ell$ with $K\leq\ell\leq \varepsilon n_Q$, we have
\begin{align*}
|\mathcal{Q}_C(n_Q,n_Q+\ell)| &\geq a_1 n_Q^{n_Q+3\ell/2-1/2}\gamma^\ell \ell^{-3\ell/2-2}\exp\left(b_1\sqrt{\ell^3n_Q^{-1}}\right);
\\
|\mathcal{Q}_C(n_Q,n_Q+\ell)|&\leq a_2 n_Q^{n_Q+3\ell/2-1/2}\gamma^\ell \ell^{-3\ell/2-2}\exp\left(b_2\sqrt{\ell^3n_Q^{-1}}\right).
\end{align*}
\end{lem}

\begin{lem}[\cite{surface}]\label{lem:complex_planar}
There exist constants  $a_3, a_4, \gamma_1, K, \varepsilon>0$ and $b_3, b_4 \in \mathbb{R}$ such that for all admissible $n_Q$ and $\ell$ with $K\leq\ell\leq \varepsilon n_Q$, we have
\begin{align*}
|\mathcal{Q}_P(n_Q,n_Q+\ell)|&\geq a_3 n_Q^{n_Q+3\ell/2-1/2}\gamma_1^\ell \ell^{-3\ell/2-3}\exp\left(b_3\sqrt{\ell^3n_Q^{-1}}\right);
\\
|\mathcal{Q}_P(n_Q,n_Q+\ell)|&\leq a_4 n_Q^{n_Q+3\ell/2-1/2}\gamma_1^\ell \ell^{-3\ell/2-3}\exp\left(b_4\sqrt{\ell^3n_Q^{-1}}\right).
\end{align*}
\end{lem}

\begin{lem}\label{lem:frac_l}
There exist constants  $c_1, c_2, K, \varepsilon>0$ and $\delta \in \left(0,1\right)$ such that for all admissible $n_Q$ and $\ell$ with $K\leq\ell\leq \varepsilon n_Q$, we have
\[
c_1^\ell \left(\frac{n_Q}{\ell}\right)^{3/2\left(\lfloor \delta \ell \rfloor- \ell\right)}\leq
\frac{|\mathcal{Q}(n_Q, n_Q+\left \lfloor \delta\ell\right \rfloor)|}{|\mathcal{Q}(n_Q, n_Q+\ell)|} \leq c_2^\ell \left(\frac{n_Q}{\ell}\right)^{3/2\left(\lfloor \delta \ell \rfloor- \ell\right)}.
\]
\end{lem}

In order to apply Lemmas \ref{lem:bound_complex} and \ref{lem:frac_l}, we need the condition $K\leq\ell\leq \varepsilon n_Q$. The next lemma shows that this is indeed not a restriction for our considerations.

\begin{lem}\label{lem:negligable}
Let $m=m(n)=n/2+s$, where $s=s(n)$ and $n^{2/3}\ll s \ll n$. Then for each $K \in \mathbb{N}$ and $\varepsilon>0$ the main contribution to $\sum_{n_Q, \ell}\tau(n_Q,\ell)$ is provided by $n_Q$ and $\ell$ with $K\leq \ell \leq \varepsilon n_Q$.
\end{lem}

In Lemma \ref{lem:unicyclic}  we observe that $u(n_U, m_U)$ stays close to one, as long as $n_U\geq 2m_U$. Thus, we will use in that case $\binom{\binom{n_U}{2}}{m_U}$ as an estimate for $\left|\mathcal{U}(n_U,m_U)\right|$. In contrast, $u(n_U, m_U)$ starts becoming quite small if $n_U<2m_U$. Hence,  in that case  we will use stronger bounds given by Lemma \ref{lem:unicyclic}\ref{lem:unicyclic_3} and \ref{lem:unicyclic_4}. Thus, we define
\[
T_1:=\sum\nolimits_{n_U\geq 2m_U}\tau(n_Q,\ell) \quad \text{and} \quad
T_2:=\sum\nolimits_{n_U< 2m_U}\tau(n_Q,\ell).
\]

\begin{lem}\label{lem:sum_1}
Let $m=m(n)=n/2+s$, where $s=s(n)$ and $n^{2/3}\ll s \ll n$. Then the main contribution to $T_1=\sum_{n_U\geq 2m_U}\tau(n_Q,\ell)$ is provided by
$n_Q=2s+O_p\left(n^{2/3}\right)$ and $\ell=\Theta\left(sn^{-2/3}\right)$.
\end{lem}

\begin{lem}\label{lem:sum_2}
Let $m=m(n)=n/2+s$, where $s=s(n)$ and $n^{2/3}\ll s \ll n$. Then the main contribution to $T_2=\sum_{n_U< 2m_U}\tau(n_Q,\ell)$ is provided by
$n_Q=2s+O_p\left(n^{2/3}\right)$ and $\ell=\Theta\left(sn^{-2/3}\right)$.
\end{lem}

Combining Lemmas \ref{lem:sum_1} and \ref{lem:sum_2} we can choose $I_{n_Q}=2s+O_p\left(n^{2/3}\right)$ and $I_\ell=\Theta\left(sn^{-2/3}\right)$. Thus, we also obtain $I_{n_C}=\Theta\left(\sqrt{n_Q\ell}\right)=\Theta\left(sn^{-1/3}\right)$. This leads to the following results on the asymptotic order of the core and excess. 

\begin{lem}\label{lem:core-excess}
Let $m=m(n)=n/2+s$, where $s=s(n)$ and $n^{2/3}\ll s \ll n$, and let $G=A(n,m)$. Then 
whp $n_C(G)=\Theta\left(sn^{-1/3}\right)$ and $ex(G)=\Theta\left(sn^{-2/3}\right)$. 
\end{lem}

In order to obtain the order of the largest component, we look at the complex part $Q_G$. Intuitively we expect that the largest component of $Q_G$ is also the largest in $G$. The following lemma tells us that this is indeed the case.

\begin{lem}\label{lem:giant_complex}
Let $m=m(n)=n/2+s$, where $s=s(n)$ and $n^{2/3}\ll s \ll n$. Moreover, let $G=A(n,m)$. Then  
$n_Q(G)-|H_1(Q_G)|=O_p\left(n^{2/3}\right).$
\end{lem}

Lemma \ref{lem:giant_complex} together with $I_{n_Q}=2s+O_p\left(n^{2/3}\right)$ implies that the complex part $Q_G$ has one component with $2s+O_p\left(n^{2/3}\right)$ vertices, while all other components are of order $O_p\left(n^{2/3}\right)$. For the non-complex components we observe that $m_U=n_U/2+O_p\left(n_U^{2/3}\right)$. Thus, for each $i\in \mathbb{N}$ the $i-$th largest non-complex component has $\Theta_p\left(n^{2/3}\right)$ vertices by Theorem \ref{thm:general} and Lemma \ref{lem:unicyclic}. This concludes the proof of Theorem \ref{thm:pt_out}.

\section{Singularity analysis: proof of Theorem \ref{thm:pt_cac}}
It suffices to show Theorem \ref{thm:pt_cac}\ref{thm:pt_cac_2}, since \ref{thm:pt_cac_1} follows from \ref{thm:pt_cac_2} and Remark 8.6 in \cite{surface}. We denote by $\mathcal{K}_c^{\circ}$ the class of connected cubic cactus weighted multigraphs, where one vertex is marked. Moreover, let $\mathcal{B}$ be the class of connected cactus weighted multigraphs, where all but one vertex have degree three and the exceptional vertex has degree two. We denote by $B(z), K(z), K_c(z)$ and $K_c^{\circ}(z)$ the exponential generating functions of the classes $\mathcal{B}, \mathcal{K}, \mathcal{K}_c$, and $\mathcal{K}_c^\circ$, respectively. By considering the marked vertex of a graph in $\mathcal{K}_c^\circ$ and distinguish some cases we obtain
\begin{align*}
K_c^\circ(z)=\frac{zB(z)}{2(1-zB(z))}+\frac{zB(z)^3}{6}.
\end{align*}
Similarly, by considering the vertex of degree two in graphs in $\mathcal{B}$ we get 
\begin{align}\label{eq:generating_function}
B(z)=\frac{z}{2(1-zB(z))}+\frac{z}{2}B(z)^2.
\end{align}
We observe that all even coefficients in $B(z)$ are zero, i.e.
$B(z)=\sum_{i\geq 1}b_{2i-1}z^{2i-1}$ for some $b_{2i-1}\in \mathbb{N}$. By taking $\widetilde{B}(u):=\sum_{i\geq 1}b_{2i-1}u^i$, we observe that 
(\ref{eq:generating_function}) translates to 
\begin{align*}
\widetilde{B}(u)=\frac{u}{2\left(1-\widetilde{B}(u)\right)}+\frac{1}{2}\widetilde{B}(u)^2.
\end{align*}
Using techniques from \cite{random_trees, anal_comb} we obtain that for $u \to r$, 
\begin{align*}
\widetilde{B}(u)=t-\rho\sqrt{1-\frac{u}{r}}+O\left(1-\frac{u}{r}\right),
\end{align*}
where $t=1-\frac{\sqrt{3}}{3}$, $r=\frac{2\sqrt{3}}{9}$, and $\rho=\frac{\sqrt{2}}{3}$. Moreover, $r$ is the unique dominant singularity of $\widetilde{B}(u)$, due to the aperiodicity of $\widetilde{B}(u)$. Next, we define $\widetilde{K}_c^\circ(u):=K_c^\circ(\sqrt{u})$,  $\widetilde{K}_c(u):=K_c(\sqrt{u})$ and $\widetilde{K}(u):=K(\sqrt{u})$. Using $u\cdot \widetilde{K}_c^\circ(u)=\widetilde{B}(u)^2-\widetilde{B}(u)^3/3$ and $K_c(z)=\int K_c^\circ(z)/z\mathrm{d}z$ we obtain that there are $k_1, k_2, k_3 \in \mathbb{R}$ such that for $u \to r$
\begin{align*}
\widetilde{K}_c(u)=k_1+k_2\left(1-\frac{u}{r}\right)+k_3\left(1-\frac{u}{r}\right)^{\frac{3}{2}}+O\left(\left(1-\frac{u}{r}\right)^2\right).
\end{align*}
Hence, there is a constant $c_1>0$ such that with $\gamma:=r^{-1/2}$ we obtain 
\[
\left[z^{2n}\right]K_c(z)=\left[u^n\right] \widetilde{K}_c(u)=c_1\gamma^{2n}n^{-\frac{5}{2}}\left(1+o(1)\right), \quad \text{as}\quad n \to \infty.
\]
Finally, we use $\widetilde{K}(u)=\exp\left(\widetilde{K}_c(u)\right)$ to obtain that there is a $c_0>0$ such that $\left[z^{2n}\right]K(z)=\left[u^n\right] \widetilde{K}(u)=c_0\gamma^{2n}n^{-\frac{5}{2}}\left(1+o(1)\right)$ for $n \to \infty$.

\section{Blocks and chords: proof of Theorem \ref{thm:out_cac}}
We will use a double counting argument to show Theorem \ref{thm:out_cac}. To that end, we need some structural information about $G=A(n,m)$. By Lemma \ref{lem:core-excess} we know that whp $n_C(G)=\Theta\left(sn^{-1/3}\right)$ and $ex(G)=\ell(G)=\Theta\left(sn^{-2/3}\right)$. Apart from that we need the two following lemmas about blocks and chords, where we call a maximal 2-connected subgraph of $G$ a \emph{block}. In addition, a \emph{chord} is an edge in $G$ that lies in a block $B$, but not in the unique Hamiltonian cycle of $B$.

\begin{lem}\label{lem:cyclic_block}
Let $m=m(n)=n/2+s$, where $s=s(n)$ and $n^{2/3}\ll s \ll n$. Then whp $A(n,m)$ does not contain a vertex that lies in three blocks.
\end{lem}

Given a chord $xy$, we denote by $B_{xy}$ the block that contains $x$ and $y$ and by $B'_{xy}$ the unique Hamiltonian cycle of $B_{xy}$. 
A chord $xy$ is said to be {\em good} (with respect to a function $h(n)=\omega(1)$) if there is a path $P_{xy}=z_0z_1 \ldots z_r z_{r+1}$ from $z_0=x$ to $z_{r+1}=y$ in $B'_{xy}$ such that 
\begin{itemize}
\item 
$z_1, \ldots, z_r$ are not endpoints of any chords in $B_{xy}$;
\item 
$r\geq n^{1/3}h(n)^{-1}+1$;
\item 
$z_i$ has degree 2 for all $i \in \mathbb{N}$ with $1\leq i \leq n^{1/3}h(n)^{-1}$.
\end{itemize}

\begin{lem}\label{lem:chord}
Let $m=m(n)=n/2+s$, where $s=s(n)$ and $n^{2/3}\ll s \ll n$ and $h(n)=\omega(1)$. Then whp $A(n,m)$ has either no chord or a good chord $xy$ (with respect to $h(n)$).
\end{lem}

Now we fix $h(n)=\omega(1)$ such that $sh(n)=o(n)$. We denote by $\mathcal{A}'(n,m)$ the subclass of $\mathcal{A}(n,m)$ containing those graphs $H$ that have a good chord, have no vertex lying in three blocks, and satisfies $n_C(H)=\Theta\left(sn^{-1/3}\right)$ and $\ell(H)=\Theta\left(sn^{-2/3}\right)$. Due to Theorem \ref{thm:pt_out} and Lemmas \ref{lem:cyclic_block} and \ref{lem:chord}, it suffices to show $|\mathcal{A}'(n,m)|=o\left(|\mathcal{A}(n,m)|\right)$. To that end, we consider the following operation for $H\in \mathcal{A}'(n,m)$: 
\begin{itemize}
\item 
We choose a good chord $xy$ and denote by $P_{xy}=z_0z_1 \ldots z_r z_{r+1}$ the corresponding good path  from $z_0=x$ to $z_{r+1}=y$.
\item
We choose $i \in \mathbb{N}$ with $1\leq i \leq n^{1/3}h(n)^{-1}$.
\item
We add the edge $z_iz_r$ and delete $z_ry$.
\end{itemize} 
We observe that we have at least $n^{1/3}h(n)^{-1}-1$ options for performing this operation. In addition, we note that the following holds in the new  graph $H'$ resulting from $H$ by the above operation: 
\begin{itemize}
\item 
$H' \in \mathcal{A}(n,m)$, $n_C(H')=n_C(H)$, and $\ell(H')=\ell(H)$;	
\item 
$z_i$ has degree 3;
\item 
$z_i$ and $z_r$ are neighbours;
\item 
there is a path from $z_i$ to $x$ such that all internal vertices have degree two;
\item 
$x$ lies in at most two blocks;
\item 
$y$ is a neighbour of $x$ such that $xy$ lies in the unique Hamiltonian cycle of the block containing $x$ and $y$.
\end{itemize}
Thus, for a fixed graph $H'$ there are at most $2\ell \cdot 3 \cdot 3 \cdot 4=\Theta\left(sn^{-2/3}\right)$ many different graphs $H$ such that we can obtain $H'$ by performing our operation in $H$. Hence, we obtain
$|\mathcal{A}'(n,m)|=O\left(\frac{sn^{-2/3}}{n^{1/3}h(n)^{-1}}\right)|\mathcal{A}(n,m)|=o\left(|\mathcal{A}(n,m)|\right)$.

\section{Sketches of proofs of auxiliary results}
\textbf{Proof of Lemma \ref{lem:estimate_core}}.
For a graph $H \in \mathcal{C}(n_C,n_C+\ell)$ we consider the following two constructions for building a graph in $\mathcal{C}(n_C+1,n_C+1+\ell)$: 
\begin{enumerate}
\item[(C1)]\label{construction_1} 
We choose an edge $e$ of $H$ which is not a chord. Then we subdivide $e$ by one vertex and label this new vertex with $n_C+1$.
\item[(C2)]\label{construction_2}
We choose a vertex $v$ in $H$ of degree 3, 4, 5 or 6 and an edge $e$ which is incident to $v$ and not a chord. Then we relabel $v$ with label $n_C+1$ and subdivide $e$ by one vertex which obtain the label of $v$.
\end{enumerate}
We observe that if $H$ has $b$ chords, then we have $n_C+\ell-b$ options for performing (C1). In addition, $H$ has at least $b/2$ vertices of degree at least three and at most $2\ell/5$ vertices of degree at least seven. Hence, we have at least $b/2-2\ell/5$ choices for performing (C2). Now if $b\leq 19\ell/20$, then we have at least $n_C+\ell/20$ choices for (C1). Otherwise if $b>19\ell/20$, then we have at least $n_C$ choices for (C1) and at least $3\ell/40$ options for (C2). We note that each graph  $H' \in \mathcal{C}(n_C+1,n_C+1+\ell)$ can be obtained at most once by performing (C1) and if this is the case, then it cannot be obtained by (C2). Finally, observing that $H'$ can be obtained at most six times by performing (C2) yields statement \ref{lem:estimate_core_1}.
\\
For \ref{lem:estimate_core_2} we call a vertex $v$ of $H' \in \mathcal{C}(n_C+1,n_C+1+\ell)$ {\em nice} if it has degree two and the two neighbours are not adjacent. We observe that $H'$ can be obtained by (C1) if the vertex $n_C+1$ is nice. We note that if $v$ has degree two and is not nice, then $v$ has a neighbour of degree at least three. Thus, $H'$ has at least $n_C+1-8\ell$ nice vertices, since the sum of all degrees of vertices of degree at least three is at most $6\ell$. As $H'$ was arbitrary, \ref{lem:estimate_core_2} follows. 
\\[0.2cm]
The statements of Lemmas \ref{lem:bound_r1}, \ref{lem:bound_r2} and \ref{lem:negligable}-\ref{lem:sum_2} are all of the type that they determine the main contribution to some sum. In order to show these results we use Lemma \ref{lem:main_analysis}, which usually requires a long and technical computation. Therefore, we provide only sketches of these proofs in this chapter, but we shall give a full proof of Lemma \ref{lem:bound_r1} in Appendix \ref{sec:appendix} to illustrate how to work out the details.
\\[0.2cm]
\textbf{Proof of Lemma \ref{lem:bound_r1} and \ref{lem:bound_r2}}.
If $\ell$ is \lq small\rq\ compared to $n_C$, we get by Lemma \ref{lem:estimate_core} that $\frac{|\mathcal{C}(n_C+1,n_C+1+\ell)|}{|\mathcal{C}(n_C,n_C+\ell)|}=n_C+\Theta(1)\ell$. Using this, we obtain $\frac{\rho(n_C+1)}{\rho(n_C)}=\left(1-\frac{n_C}{n_Q}\right)\left(1+\Theta(1)\frac{\ell}{n_C}\right)$. Hence, we expect that the main contribution to (\ref{eq:2}) is provided by terms with $n_C=\Theta\left(\sqrt{n_Q\ell}\right)$.
\\[0.2cm]
\textbf{Proof of Lemma \ref{lem:bound_complex}}.
Combining Lemmas \ref{lem:bound_r1} and \ref{lem:bound_r2} together with (\ref{eq:frac_r}) we obtain
\begin{align*}
\frac{|\mathcal{Q}(n_Q+1, n_Q+1+\ell)|}{|\mathcal{Q}(n_Q, n_Q+\ell)|}&\approx \frac{\rho_1\left(\sqrt{n_Q\ell}\right)}{\rho_0\left(\sqrt{n_Q\ell}\right)}\\
&\approx (n_Q+1)\exp\left(\frac{\sqrt{n_Q\ell}}{n_Q-\sqrt{n_Q\ell}+1}+\frac{n_Q-\sqrt{n_Q\ell}-1}{n_Q}\right)
\\
&\approx
(n_Q+1)\exp\left(1+\frac{\ell}{n_Q}\right).
\end{align*}
\textbf{Proof of Lemma \ref{lem:complex_cacti}}.
Using the core-kernel approach from \cite{surface} and following the lines of the proofs of Lemma 4.9(ii) and Corollary 4.11 in \cite{surface} yields the assertion. (A detailed proof can be found in Appendix \ref{sec:appendix2}).
\\[0.2cm]
\textbf{Proof of Lemma \ref{lem:frac_l}}.
We note that $\exp\left(\Theta(1)\sqrt{\ell^3n_Q^{-1}}\right)=\exp\left(\Theta\left(1\right)\ell\right)$, since $\ell=O\left(n_Q\right)$. Then the statement follows by combining Lemmas \ref{lem:complex_cacti} and \ref{lem:complex_planar} together with (\ref{eq:bound_complex}).
\\[0.2cm]
\textbf{Proof of Lemma \ref{lem:negligable}}. 
We denote by $\mathcal{T}$ the class of cactus graphs. Clearly, we have $\left|\mathcal{A}(n,m)\right|\geq\left|\mathcal{T}(n,m)\right|$, because every cactus graph is also an outerplanar graph. By the core-kernel approach we obtain that there is a $c>0$ and $N\in \mathbb{N}$ such that
$
\left|\mathcal{T}(n,m)\right|\geq \frac{n^{n-1/2}}{\left(n-2s\right)^{n/2-s}}\exp\left(\frac{n}{2}-s+c \cdot \frac{s}{n^{2/3}}\right) 
$ for all $n\geq N$.
On the other hand, we can bound $\tau(n_Q,\ell)$ by Lemmas \ref{lem:unicyclic} and \ref{lem:complex_planar}. By doing so we obtain that 
$\sum_{\ell<K, n_Q}\tau(n_Q,\ell)=o\left(\left|\mathcal{T}(n,m)\right|\right)$.
Hence, the terms provided by $\ell<K$ are negligible in $\sum_{n_Q, \ell}\tau(n_Q,\ell)$. Similarly, one can also show that this is true for the terms provided by $\ell>\varepsilon n_Q$.
\\[0.2cm]
\textbf{Proof of Lemma \ref{lem:sum_1}}.
By Lemma \ref{lem:unicyclic} we may consider $Y_1=\sum_{n_U\geq 2m_U}\upsilon_1(n_Q,\ell)$ instead of $T_1$, where $\upsilon_1(n_Q,\ell):=\binom{n}{n_Q}\left|\mathcal{Q}(n_Q,n_Q+\ell)\right|\binom{\binom{n_U}{2}}{m_U}$. Then we obtain by Lemma \ref{lem:frac_l} that 
$\frac{\upsilon_1\left(n_Q,\left \lfloor \delta\ell\right \rfloor\right)}{\upsilon_1(n_Q,\ell)}=\Big(\Theta(1)\frac{n_Q^{3/2}m_U}{\ell^{3/2}n^2}\Big)^{\lfloor \delta \ell \rfloor- \ell}.$ 
Thus, the main contribution to $Y_1$ is provided by $n_Q$ and $\ell$ with $\frac{n_Q}{\ell}=\Theta\big(n^{4/3}m_U^{-2/3}\big)$. Combining that together with Lemma \ref{lem:bound_complex} we get
\begin{align*}
\frac{\upsilon_1(n_Q+1,\ell)}{\upsilon_1(n_Q,\ell)}=\exp\left(O\left(n^{-2/3}\right)-\Theta(1)\left(1-\frac{2m_U}{n_U}\right)^2\right).
\end{align*}
Thus, the main contribution to $Y_1$ is provided by $n_Q$ and $\ell$ with $\frac{n_Q+2\ell-2s}{n_U}=\left(1-\frac{2m_U}{n_U}\right)=O_p\left(n^{-1/3}\right)$, which yields $n_Q+2\ell-2s=O_p\left(n^{2/3}\right)$. Together with $\frac{n_Q}{\ell}=\Theta\big(n^{4/3}m_U^{-2/3}\big)$ this implies $n_Q=2s+O_p\left(n^{2/3}\right)$ and $\ell=\Theta\left(sn^{-2/3}\right)$.
\\[0.2cm]
\textbf{Proof of Lemma \ref{lem:sum_2}}.
We define
\begin{align*}
\upsilon_2(n_Q,\ell):=\upsilon_1(n_Q,\ell)c\left(\frac{2}{e}\right)^{2m_U-n_U}\frac{m_U^{m_U+1/2}n_U^{n_U-2m_U+g(n_Q)}}{\left(n_U-m_U\right)^{n_U-m_U+1/2}},
\end{align*}
where $c>0$, $h(n)=\omega(1)$, $\upsilon_1(n_Q,\ell)$ as in the proof of Lemma \ref{lem:sum_1} and $g(n_Q):=\frac{1}{2}$ if $n_Q\leq 2s-n^{2/3}h(n)$ and $g(n_Q):=0$ otherwise. By Lemma \ref{lem:unicyclic} we can choose $h(n)$ and $c$ so that for all admissible $n_Q$ and $\ell$, we have $\tau(n_Q,\ell)\leq \upsilon_2(n_Q,\ell)$. Similarly as in the proof of Lemma \ref{lem:sum_1} we obtain that the main contribution to $Y_2:=\sum_{n_U<2m_U}\upsilon_2(n_Q,\ell)$ is provided by $n_Q=2s+O_p\left(n^{2/3}\right)$ and $\ell=\Theta\left(sn^{-2/3}\right)$. For such $n_Q$ and $\ell$ we have $g(n_Q)=0$ and by Lemma \ref{lem:unicyclic}\ref{lem:unicyclic_2} $|\mathcal{U}(n_U,m_U)|=\Theta_p(1)\binom{\binom{n_U}{2}}{m_U}$. Using that we obtain $\frac{\upsilon_2(n_Q,\ell)}{\tau(n_Q,\ell)}=\Theta_p(1)$, which shows the statement.
\\[0.2cm]
\textbf{Proof of Lemma \ref{lem:giant_complex}}.
Let $\widetilde{n}=n_Q-|H_1(Q_G)|$ and we look at the following operation in $G$. We add an edge between two different complex components and delete an edge in a non-complex component. We have whp $\Omega\left(s\widetilde{n}n\right)$ choices for performing this operation. We observe that in the reverse operation we delete an edge from the core and add some edge. We can do that whp in $O\left(sn^{-1/3}n^2\right)$ different ways. Hence, it follows that $\widetilde{n}=O_p\left(\frac{sn^{5/3}}{sn}\right)=O_p\left(n^{2/3}\right)$.
\\[0.2cm]
\textbf{Proof of Lemma \ref{lem:cyclic_block}}.
Let $H \in \mathcal{A}(n,m)$ be a graph that has a vertex lying in three blocks. We consider the following operation in the core $C_H$: 
\begin{itemize}
\item 
We choose a vertex $x$ that lies in three blocks;
\item 
Let $X$ be the component of $C_H$ containing $x$. Then we choose a component $Y$ of $X-x$ that contains at most $n_C(H)/3$ vertices, but two neighbours of $x$ (in $H$);
\item 
We choose a vertex $y$ in $C_H$ which is not in $Y$ and has degree two;
\item 
For all neighbours $z$ of $x$ in $Y$ we delete the edge $xz$ and insert the edge $yz$.
\end{itemize}
We observe that we have at least $2n_C(H)/3-2\ell=\Theta\left(sn^{-1/3}\right)$ options for performing this operation. On the other hand, we note that in a constructed graph $H'$ the following holds:
\begin{itemize}
\item 
$H' \in \mathcal{A}(n,m)$, $n_C(H')=n_C(H)$ and $\ell(H')=\ell(H)$;		
\item
$y$ lies in one or two blocks and has at least degree four;
\item 
$x$ has at least degree four.
\end{itemize}
Hence, a fixed graph $H'$ can be constructed in at most $2\ell \cdot 2 \cdot 2\ell=\Theta\left(s^2n^{-4/3}\right)$ many different ways. Now the statement follows, since $\Theta\left(\frac{s^2n^{-4/3}}{sn^{-1/3}}\right)=o(1)$.
\\[0.2cm]
\textbf{Proof of Lemma \ref{lem:chord}}.
We consider the kernel $K_H$ of a graph $H\in \mathcal{A}(n,m)$ which has a chord. Then $K_H$ has a chord $xy$ with the following property: If $B'$ is the unique Hamiltonian cycle of the block $B$ containing $x$ and $y$, then there is a path $z_0=x, z_1, \ldots, z_{t+1}=y$ in $B'$ such that there is no chord in $B$ containing one of the vertices $z_1, \ldots, z_t$. Next, we choose a random core which can be obtained by subdividing the edges of $K_H$ which are not chords by $n_C(H)-|K_H|$ additional vertices. We denote by $X$ the number of vertices which subdivide the edge $z_0z_1$. Using a \lq bins and balls\rq\ type argument, we can show that  $\mathbb{P}\left[X=j\right]\leq \mathbb{P}\left[X=0\right]$ for any $j \in \mathbb{N}$ and $\mathbb{P}\left[X=0\right]=O\left(\frac{\left|K_H\right|}{n_C(H)-\left|K_H\right|}\right)=O\left(n^{-1/3}\right)$. Thus, $\mathbb{P}\left[X\leq n^{1/3}h(n)^{-1}\right]\leq \left(n^{1/3}h(n)^{-1}+1\right)\mathbb{P}\left[X=0\right]=o(1)$, i.e. whp  $z_0z_1$ is subdivided by at least $n^{1/3}h(n)^{-1}+1$ vertices, which shows the statement.

\section{Concluding remarks}
Kang, Mo{\ss}hammer, and Spr\"ussel \cite{surface} showed that graphs on orientable surfaces feature a second phase transition at $m=n+O\left(n^{3/5}\right)$, where the number of vertices outside the largest component becomes sublinear. By Theorem \ref{thm:pt_cac} and Remark 8.6 in \cite{surface} this is also true for random cactus graphs. Thus, we believe that this should also be the case for random outerplanar graphs, since the class of outerplanar graphs lies \lq between\rq\ the class of cactus graphs and the class of graphs on orientable surfaces. Unfortunately, our method does not seem to work when $m=n+o(n)$. This is mainly because the bound in Lemma \ref{lem:frac_l} is not good enough in that regime.

Theorem \ref{thm:out_cac} raises the following question. How does the probability that $A(n,m)$ is a cactus graph behave if $m$ grows? By looking at the proof of Theorem \ref{thm:out_cac} a natural guess would be the following.

\begin{conjecture}
If $m=\alpha n$ for $1/2<\alpha<1$, then the probability that $A(n,m)$ is a cactus graph is bounded away from 0 and 1.
\end{conjecture}

\begin{conjecture}
If $m=n+t$ for $t=o(n)$, then whp $A(n,m)$ is not a cactus graph.
\end{conjecture}

\section*{Acknowledgement}
This extended abstract is based on the Master's thesis of Michael Missethan \cite{master_thesis}. Detailed proofs of main and auxiliary results of this paper can be found in \cite{master_thesis}.

\bibliographystyle{plain}

\appendix
\newpage
\section{An application of Lemma \ref{lem:main_analysis}: proof of Lemma \ref{lem:bound_r1}}\label{sec:appendix}
To illustrate how to apply Lemma \ref{lem:main_analysis} we prove Lemma \ref{lem:bound_r1} in this section (the proof of Lemma \ref{lem:bound_r2} is similar). We start by getting an upper bound for $\frac{\rho\left(n_C\right)}{\rho\left(n_C+1\right)}$. By Lemma \ref{lem:estimate_core}\ref{lem:estimate_core_1} we obtain
\begin{align*}
\frac{\rho\left(n_C\right)}{\rho(n_C+1)}&=\frac{n_C+1}{n_Q-n_C}\cdot \frac{n_Cn_Q}{n_C+1} \cdot \frac{|\mathcal{C}(n_C, n_C+\ell)|} {|\mathcal{C}(n_C+1, n_C+1+\ell)|}\\&\leq
\frac{n_Cn_Q}{n_Q-n_C} \frac{1}{n_C+\frac{\ell}{80}}
\\
&=\left(1- \frac{\ell}{80n_C+\ell}\right)\left(1+\frac{n_C}{n_Q-n_C}\right)
\\
&\leq \exp\left(-\frac{\ell}{80n_C+\ell}+\frac{n_C}{n_Q-n_C}\right).
\end{align*}
Next, we observe that $\ell\leq n_C\leq n_Q$, since an outerplanar graph on $n_C$ vertices can have at most $2n_C$ edges. Hence, we can choose $c>0$ small enough such that for all $n_C\leq 2c\sqrt{n_Q\ell}$
\begin{align*}
\frac{\rho\left(n_C\right)}{\rho(n_C+1)}&\leq \exp\left(-\frac{\ell}{81n_C}+\frac{2n_C}{n_Q}\right)
\\
&\leq \exp\left(-\frac{\ell}{81\cdot 2c\sqrt{n_Q\ell}}+\frac{2\cdot 2c\sqrt{n_Q\ell}}{n_Q}\right)
\\
&\leq\exp\left(-\sqrt{\frac{\ell}{n_Q}}\right)=\exp\left(h(n)\right),
\end{align*}
where $h(n):=-\sqrt{\frac{\ell}{n_Q}}$. We also define $g(n):=c\sqrt{n_Q\ell}$ and $f_n(n_C):=n_C+g(n)$. Then we obtain for all $n_C\leq c\sqrt{n_Q\ell}$
\begin{align*}
\frac{\rho\left(n_C\right)}{\rho\left(f_n(n_C)\right)}=\prod_{k=n_C}^{f_n(n_C)-1}\frac{\rho(k)}{\rho(k+1)}\leq \exp\left(h(n)\right)^{g(n)}=\exp\left(-c\ell\right).
\end{align*}
Finally, that yields
\begin{align*}
\sum_{n_C\leq c\sqrt{n_Q\ell}}\rho(n_C)\leq \exp\left(-c\ell\right) \sum_{n_C\leq c\sqrt{n_Q\ell}}\rho(f_n(n_C))\leq \exp\left(-c\ell\right) \sum_{n_C}\rho(n_C),
\end{align*}
which shows the statement.
\\
We conclude this section by observing an immediate consequence of Lemma \ref{lem:bound_r1}. Assuming $\ell =\omega(1)$, which is true due to Lemmas \ref{lem:sum_1} and \ref{lem:sum_2}, we have $\frac{\rho\left(n_C\right)}{\rho\left(f_n(n_C)\right)}\leq\exp\left(-c\ell\right)=o(1)$. Then Lemma \ref{lem:main_analysis} implies that the terms provided by $I_1(n):=\left\{n_C~|~n_C\leq c\sqrt{n_Q\ell}\right\}$ are negligible in $\sum_{n_C}\rho(n_C)$.

\section{Proof of Lemma \ref{lem:complex_cacti}}\label{sec:appendix2}
We shall focus on the proof of the lower bound, since the upper bound can be shown in a similar way. We will use the core-kernel approach from \cite{surface} and recall that $\mathcal{T}$ is the class of all cactus graphs. Then we denote by $\mathcal{C}_C$ the class of all cores of graphs in $\mathcal{T}$ and by $\mathcal{K}_C$ the class of all kernels of graphs in $\mathcal{T}$. Analogously to (\ref{eq:2}) we obtain
\begin{align}\label{eq:complex_cacti}
|\mathcal{Q}_C(n_Q,n_Q+\ell)|=\sum\nolimits_{n_C}\binom{n_Q}{n_C}|\mathcal{C}_C(n_C, n_C+\ell)|n_Cn_Q^{n_Q-n_C-1}.
\end{align}
We claim that 
\begin{align}\label{ineq:complex_cacti}
|\mathcal{C}_C(n_C, n_C+\ell)|\geq \binom{n_C}{2\ell}\left|\mathcal{K}_C(2\ell, 3\ell)\right|\left(n_C-2\ell\right)!\binom{n_C-5\ell-1}{3\ell-1}.
\end{align}
Indeed, we can construct (not necessarily all) graphs from $\mathcal{C}_C(n_C, n_C+\ell)$ in the following way. We choose $2\ell$ labels from $[n_C]$ for the vertices of the kernel. Then we pick a kernel $K$ from $\mathcal{K}_C(2\ell, 3\ell)$ and assign the labels chosen before to the vertices of $K$. Finally, we subdivide the edges of the kernel by the $(n_C-2\ell)$ remaining vertices such that each edge is subdivided by at least two vertices, which guarantees that the obtained graph is simple. Thus, all constructed graphs are in $\mathcal{C}_C(n_C, n_C+\ell)$. We note that there are $w(K)\left(n_C-2\ell\right)!\binom{n_C-5\ell-1}{3\ell-1}$ many ways to get such a subdivision, where $w(K)=2^{-e_1(K)-e_2(K)}$ and $e_1(K)$ denotes the number of loops in $K$ and $e_2(K)$ the number of double edges in $K$. In addition, we note that in $\left|\mathcal{K}_C(2\ell, 3\ell)\right|$ each kernel $K$ is counted with a weight of $w(K)$. Then inequality (\ref{ineq:complex_cacti}) follows by the aforementioned construction. Combining (\ref{eq:complex_cacti}) and (\ref{ineq:complex_cacti}) we obtain
\begin{align}\label{eq:complex_cacti2}
	|\mathcal{Q}_C(n_Q,n_Q+\ell)|&\geq\frac{\left|\mathcal{K}_C(2\ell, 3\ell)\right|n_Q^{n_Q-1}}{(2\ell)!(3\ell-1)!}\sum\nolimits_{n_C}\left(n_Q\right)_{n_C}\left(n_C-5\ell-1\right)_{3\ell-1}n_Cn_Q^{-n_C}\nonumber
	\\
	&=\frac{\left|\mathcal{K}_C(2\ell, 3\ell)\right|n_Q^{n_Q-1}}{(2\ell)!(3\ell-1)!}\sum\nolimits_{n_C}\nu(n_C),
\end{align}
where $\nu(n_C):=\left(n_Q\right)_{n_C}\left(n_C-5\ell-1\right)_{3\ell-1}n_Cn_Q^{-n_C}$. Next, we observe that
\begin{align*}
\frac{\nu(n_C+1)}{\nu(n_C)}=\frac{n_Q-n_C}{n_Q}\frac{n_C-5\ell}{n_C-8\ell+1}\frac{n_C+1}{n_C}.
\end{align*}
We note that $\frac{\nu(n_C+1)}{\nu(n_C)}$ is decreasing in $n_C$ and that $\frac{\nu(\overline{n_C}+1)}{\nu(\overline{n_C})}\approx 1$ for $\overline{n_C}=\sqrt{3n_Q\ell}$. Thus, we expect that we obtain a good approximation for $\sum\nolimits_{n_C}\nu(n_C)$ by considering only terms whose index is \lq close\rq\ to $\overline{n_C}$. In the following we make that more precise. We note that for $\ell\leq\varepsilon n_Q$ and $\varepsilon>0$ small enough, we get
\begin{align}\label{eq:complex_cacti3}
\nu\left(\overline{n_C}\right)\geq \left(1-\frac{\sqrt{3n_Q\ell}}{n_Q}\right)^{\sqrt{3n_Q\ell}}\left(\sqrt{3n_Q\ell}-8\ell\right)^{3\ell}\geq \exp\left(-6\ell\right)\sqrt{n_Q\ell}^{3\ell}.
\end{align}
Next, we distinguish two cases. First we assume $\ell\leq \sqrt{n_Q}$. Then we get for all $n_C\geq\overline{n_C}-\sqrt{n_Q}$ and $\varepsilon>0$ small enough
\begin{align*}
\frac{\nu({n_C}+1)}{\nu\left({n_C}\right)}&\leq\left(1-\frac{{n_C}}{n_Q}\right)\left(1+\frac{3\ell}{{n_C}-8\ell}\right)\left(1+\frac{1}{{n_C}}\right)
\\
&\leq \exp\left(-\sqrt{\frac{3\ell}{n_Q}}+\frac{3\ell}{\sqrt{3n_Q\ell}-\sqrt{n_Q}-8\ell}+\frac{3}{\sqrt{n_Q}}\right)
\\
&\leq \exp\left(\sqrt{\frac{3\ell}{n_Q}}\cdot\frac{27\sqrt{n_Q}}{\sqrt{3n_Q\ell}}+\frac{3}{\sqrt{n_Q}}\right)= \exp\left(\frac{30}{\sqrt{n_Q}}\right).
\end{align*}
Hence, we obtain $\nu\left(n_C\right)\geq\nu\left(\overline{n_C}\right)\exp(-30)$ for all $\overline{n_C}-\sqrt{n_Q}\leq n_C\leq\overline{n_C}$. Combining that together with (\ref{eq:complex_cacti2}), (\ref{eq:complex_cacti3}) and Theorem \ref{thm:pt_cac} yields 
\begin{align*}
	|\mathcal{Q}_C(n_Q,n_Q+\ell)|
	&\geq\frac{\left|\mathcal{K}_C(2\ell, 3\ell)\right|n_Q^{n_Q-1}}{(2\ell)!(3\ell-1)!}\sqrt{n_Q}\nu\left(\overline{n_C}\right)\exp(-30)
	\\
	&\geq \Theta(1)^{\ell}n_Q^{n_Q+3\ell/2-1/2}\ell^{3\ell/2-5/2-3\ell+1/2}
	\\
	&=\Theta(1)^{\ell}n_Q^{n_Q+3\ell/2-1/2}\ell^{-3\ell/2-2},
\end{align*}
which shows the statement for the case $\ell\leq \sqrt{n_Q}$. Finally, we assume $\ell>\sqrt{n_Q}$. Then we get by (\ref{eq:complex_cacti2}), (\ref{eq:complex_cacti3}) and Theorem \ref{thm:pt_cac} for $\varepsilon>0$ small enough
\begin{align*}
	|\mathcal{Q}_C(n_Q,n_Q+\ell)|
	&\geq\frac{\left|\mathcal{K}_C(2\ell, 3\ell)\right|n_Q^{n_Q-1}}{(2\ell)!(3\ell-1)!}\nu\left(\overline{n_C}\right)
	\\
	&\geq
	\Theta(1)^\ell n_Q^{n_Q-1+3\ell/2}\ell^{3\ell/2-5/2-3\ell+1/2}
	\\
	&=\Theta(1)^\ell n_Q^{n_Q+3\ell/2-1/2}\ell^{-3\ell/2-2}n_Q^{-1/2}
	\\
	&\geq
    \Theta(1)^\ell n_Q^{n_Q+3\ell/2-1/2}\ell^{-3\ell/2-2}\exp\left(-\sqrt{\frac{\ell^3}{n_Q}}\right),
\end{align*}
as desired.
\end{document}